# DISCUSSION OF "BREAKDOWN AND GROUPS" BY P. L. DAVIES AND U. GATHER

By Frank Hampel

*ETH Zürich*

**1. Introductory remarks.** It is a great pleasure for me to be invited to comment upon the nice and elegant and in parts thought-provoking paper by Davies and Gather. The authors also asked me specifically to comment upon the historical roots of the breakdown point (BP), and my thoughts about it. I shall try to do so, stressing in particular aspects and work that are not published.

**2. Some thoughts with the definition of the breakdown point.** In my thesis [Hampel (1968)] I developed what was later also called the "infinitesimal approach to robustness," based on one-step Taylor expansions of statistics viewed as functionals (the "influence curves" or "influence functions"), a technology which for ordinary functions has long been indispensable in engineering and the physical sciences, and also for much theoretical work. However, it was always clear to me that this technology needed to be supplemented by an indication up to what distance (from the model distribution around which the expansion takes place) the linear expansions would be numerically, or at least semiquantitatively, useful. The simplest idea that came to my mind (simplicity being a virtue, also in view of Ockham's razor) was the distance of the nearest pole of the functional (if it was unbounded); see the graphs in Hampel, Ronchetti, Rousseeuw and Stahel [(1986), pages 42, 48, 177]. Thus, right after defining the "bias function" (without using this term) as the (more complicated) bridge between model and pole, I introduced the "break-down point" on page 27 (Chapter C.4) of my thesis and, in a slight variant (by not requiring qualitative robustness anymore and therefore treating it as a purely global concept), as "breakdown point" on page 1894 in Hampel (1971). I was, of course, clearly inspired by Hodges (1967), whose intuition went in a similar direction, and by his "tolerance of extreme values"; however, his concept is not only much more limited, it is formally not even a special case of the breakdown point. [And contrary to

---

a claim someone spread later, the term "breakdown point" does not occur anywhere in Hodges (1967).]

My definition of the BP is asymptotic, because I believe that suitable, elegant and properly interpreted asymptotics is much more informative and more generally applicable than specific or even clumsy finite-sample definitions. However, I also believe that asymptotic results need interpretations (and often numerical checks) in finite-sample frameworks, and lack of this may even be the biggest gap separating mathematical statistics from good applications of statistics [cf. Hampel (1998)].

Since I consider the finite-sample interpretations of an asymptotic definition (even different ones under different circumstances) an integral part of the properly interpreted definition, I never felt the need to introduce a general explicit definition of a finite-sample breakdown point. In fact, different needs require different definitions. In my eyes, the BP should be a flexible tool adapted to the requirement of specific problems (see also below).

Informal finite-sample BPs have been used in Andrews et al. (1972), and, for example, in many of my papers, starting with Hampel (1973). Often, the lower (or else upper) gross-error finite-sample BP is sufficient. But Grize (1978) showed the need for the "total-variation BP" in a specific situation concerning correlations. A standard reference is Donoho and Huber (1983). But in the background remains the fact that the BP is originally defined with the Prohorov distance. Very often we can forget this somewhat awkward distance and simplify; but whenever it is needed, we have to be ready to dig it out again.

The use of the Prohorov distance needs some explanation, also in view of the paper under discussion. Many good mathematical statisticians strive for the greatest generality, without regarding the practical implications. In some way, this is legitimate (and even required by the mathematical side of statistics). But I rather try to find the specific concepts most suitable for the problem at hand. Thus, as I explained elsewhere [e.g., in Hampel (1968)], I find it necessary to use the weak (formerly weak*) topology for general robustness problems, which is metrized by the Prohorov (former spelling Prokhorov) distance, which in turn has a nice interpretation in terms of the model deviations occurring in real life. [For more technical details, see Huber (1981).] This does not preclude the possibility of simplifying in specific situations. For example, Huber's (1964) gross-error model is theoretically too narrow in scope, but it captures the most important deviations from the idealized model, and the solutions found are also useful and good in the more general situation [cf. Hampel (1992)].

Note that there is nothing about equivariance in my definition. If it is available, it simplifies life often tremendously and allows a nice mathematical theory with beautiful theorems; but I do not consider it an intrinsically necessary part of a general statistical theory (cf. also Fisher's view on his general theory of estimation).

It may be considered ironic in view of the present discussion paper that in my definition of the BP (with the compact proper subset of the parameter space), I explicitly thought of correlation statistics as examples, where there is no equivariance at all. As the authors correctly observe, such BPs have not become popular at all (so far), giving some credit to their stress on group structures. Compare also below.

**3. Some further developments.** The above idea of combining linear extrapolation with the BP was very successful in the cases tried [cf. Hampel et al. (1986), Subsection 1.3e, in particular Table 1 on page 50, which reproduces Huber's (1964) Table I—and thus his minimax results—very accurately; and Hampel (1983), page 214, which reproduces some of the Monte Carlo results in Andrews et al. (1972)]. As a rule of thumb, under mild conditions the linear extrapolations seem to be very accurate up to $BP/4$, and still numerically quite usable even somewhat beyond $BP/2$.

Another for me quite surprising success was the explanation of the (partly unsuspectedly bad) empirical behavior of various rejection rules just by means of the BP [see Hampel (1985)].

For regression I introduced the conditional BP given the design in Hampel [(1975), page 379] (implicitly and condensed because of the page limit imposed). It is more intricate, but also more informative (once the design is fixed or the data are in) than the unconditional BP [which was mostly used later on, except, e.g., in Hampel et al. (1986), page 328, unfortunately without stressing the difference between the two BPs].

Some definitions of variants of the BP, adapted to specific ANOVA-type problems, have been given by Hampel (1987), by Mili, Phaniraj and Rousseeuw (1990) and by Ruckstuhl (1995); see also Stahel, Ruckstuhl, Senn and Dressler (1994).

The BP seemed to be trivial, with $BP = 50\%$ easily possible in the models considered, until Maronna (1976) essentially showed the upper bound to be $= 1/\text{dimension}$ for "nice" equivariant estimators in multivariate and multiple regression situations. Much effort has since been put into keeping the equivariance and reaching $BP = 50\%$ with "pathological" estimators [the first prototype having been the "shordth" or "minimum median deviation" method in Hampel (1975), page 380, later popularized under the name "least median of squares"]. But from a practical point of view, I find it more reasonable to give up exact equivariance. Gross errors are often partly in single coordinates and are not equivariant, even if the ideal model is.

For general nonlinear models, equivariance may not be attainable at all, but it may make perfect sense to look at the ("a posteriori") BP at and in a neighborhood of the fitted model. Compare also the highly condensed first sentence of 3.3 on page 380 in Hampel (1975), valid also under nonequivariance (and suggesting nice quantitative theorems under equivariance).

**4. The thesis by Grize.** In his unpublished Diplomarbeit, Yves-Laurent Grize (1978) made a thorough and deep investigation of the breakdown properties and influence functions of correlation measures, notably of the Kendall (K), Spearman (S) and quadrant (Q) rank correlations. He noted that the BP actually depends on the model F, and that also the specification of the "distance" may make a difference. For some F's, $BP(K) = BP(S) = 1$, while for others $BP(K) = (3/2)BP(S) < 1$, and in again another situation $BP(K) = 0.29$, $BP(Q) = 0.25$ and $BP(S) = 0.21$. Grize showed that for correlations, the gross-error BP is often not suitable, and he used the better-fitting total-variation BP instead. He briefly also discussed the possibility of the (much more complicated) Prohorov-distance BP, and of ranks (by gross mistakes) outside the range from 1 to $n$. It appears that often Kendall's rank correlation is considerably more robust than Spearman's (and that there are some meaningful numbers and results to be taken out for statistical practice), but a lot depends on the precise specification of the situation.

My first reaction was disappointment. The results were just not as simple and beautiful as we then were used to in robust statistics. But the thesis is a valuable piece of work, and I regret very much that by some unfortunate circumstances it never found its way into the printed literature. Perhaps the time was not yet ripe for it. It seems that in recent years, some fragments of it are being rediscovered [cf. Bin Abdullah (1990) and Dehon and Croux (2003)], partly with seemingly contradictory results ("BP small" vs. "BP = 1"), which may be due to insufficient care for the fine details (which really matter here). As the recent interest in the (formerly "too complicated") "bias function" (cf. above) shows, it may well be that in the near future "complicated" BPs without a natural equivariance structure will become more popular.

**5. The small print.** It seems to me that in the regression-through-0 example of the discussion paper, there is the same play with asymptotics (concerning both BP and consistency) going on which Fisher [(1956); cf. also Hampel (1998)] complained about when he defended his definition of consistency against Neyman's. In the case of the two location samples, I guess that the unnamed estimator not breaking down under the specific model and alternatives of (6.2) has a BP of 0.005—if the Prohorov distance is taken into account, which in such situations cannot be neglected. Thus, it really seems to be a case of small print that has been forgotten.

**Acknowledgment.** The author is grateful to U. Gather for digging out some references.

DISCUSSION OF "BREAKDOWN AND GROUPS" BY P. L. DAVIES AND U. GATHER 5

Seminar für Statistik
ETH Zürich
CH-8092 Zürich
Switzerland
e-mail: hampel@stat.math.ethz.ch